\documentclass[12pt,leqno]{article}
\input xy
\xyoption{all}
\usepackage{theorem,amssymb,amsmath}

\begin{document}
\title{Complex cobordism and embeddability of CR-manifolds}
\author{Bruno De Oliveira \thanks {Partially supported by a NSF Postdoctoral Research Fellowship}} 
\maketitle

\

This paper studies complex cobordisms between compact, three dimensional, strictly pseudoconvex Cauchy-Riemann manifolds. Suppose the complex cobordism is given by a complex 2-manifold $X$ with one pseudoconvex and one pseudoconcave end. We answer the following questions. To what extent is $X$ determined by the pseudoconvex end? What is the relation between the embeddability of the pseudoconvex end  and the embeddability of the pseudoconcave end of $X$? Do all CR-functions on the pseudoconvex end of $X$ extend to holomorphic functions on the interior of $X$? \par
We introduce two methods to construct pseudoconcave surfaces that show that the complex 2-manifold $X$ giving a complex cobordism is not determined by the pseudoconvex end. These two constructions give new methods to construct non-embeddable Cauchy-Riemann 3-manifolds and prove that embeddability of a strictly pseudoconvex Cauchy-Riemann 3-manifold is not a complex-cobordism invariant. We show that a new phenomenon occurs: there are CR-functions on the pseudoconvex end that do not extend to holomorphic functions on $X$. We also show that the extendability of the CR-functions from the pseudoconvex end is necessary but not sufficient for embeddability to be preserved under complex cobordisms.\par

A compact (2n+1)-dimensional Cauchy-Riemann manifold (CR-(2n+1)-manifold) consists of: a compact (2n+1)-dimensional manifold $M$, a rank n complex subbundle $T^{0,1}M \subset TM \otimes \mathbb C$ satisfying $T^{1,0}M \cap T^{0,1}M=\{0\}$ ($T^{1,0}M \equiv \overline {T^{0,1}M}$) and the integrability condition $[ \overline Z,\overline {Z'}] \in T^{0,1}M$ for local sections $\overline Z,\overline {Z'} \in {T^{0,1}M}$. If the CR-(2n+1)-manifold has the additional property that any nonvanishing local section $\overline {Z}$ of $T^{0,1}M$ is such that $[Z,\overline {Z}] \not \in T^{0,1}M \oplus T^{1,0}M$, it is called strictly pseudoconvex (SPCR-(2n+1)-manifold). \par

A differentiable function $f:M \to \mathbb C$ is said to be a CR function if it verifies the tangential Cauchy-Riemann equations $\overline Z f=0$ for all local sections of $T^{0,1}M$. In other words, CR functions are the elements  of the kernel of the $\overline {\partial}_b$ operator defined as $\overline {\partial}_b f=df|_{T^{0,1}M}$.\par

\

{\bf Definition} A CR-(2n+1)-manifold $M$ is embeddable if there is a collection of CR functions $\{f_1,...,f_N\}$ giving an embedding $f=(f_1,...,f_N):M \to \mathbb C^N$.

\

The origin of CR-structures lies on the study of real hypersurfaces in $\mathbb C^N$.
The CR-structure being given by $T^{0,1}M=T^{0,1}\mathbb C^N \cap TM \otimes \mathbb C$.
A CR-function is a function that satisfies the differential conditions of holomorphy that can be verified along M. An  important example of SPCR-manifolds are the boundaries of strictly pseudoconvex domains in $\mathbb C^N$.

\

{\bf Definition} A CR-(2n+1)-manifold $M$ is fillable if it is the boundary of a normal complex Stein space $X$ with the CR-structure induced from the complex structure of $X$. A filling of $M$ is any normal complex space $Y$ which is a modification of a Stein space and whose boundary is $M$.

\

An example of Rossi [Ro65] showed that not all compact  SPCR-3-manifolds are embeddable. Later Boutet de Monvel [Bou74] showed that for $n\ge 2$ any compact SPCR-(2n+1)-manifold $M$ is embeddable, this follows from the result of Kohn [Ko86] stating that the range of $\overline {\partial}_b$ is closed plus the result of  [Bou74] stating that if the range of $\overline {\partial}_b$ is closed then $M$ is embeddable. This same result implies that if $M$ is fillable then $M$ is embeddable since in [KoRo65] it was shown that if a SPCR-manifold  is fillable then the range of $\overline {\partial}_b$ closed. In the other direction, it follows from the results of Harvey-Lawson [HaLa75]  that every embeddable compact SPCR-(2n+1)-manifold is fillable by a unique normal Stein space. A significant amount of work has been produced to understand the structure of the set of embeddable SPCR-structures for 3-manifolds by Bland, Burns, Catlin, Epstein, Henkin and Lempert. \par

The results described above state that there is a rigidity on the possible fillings of an embeddable SPCR-3-manifold $M_0$, since any filling of $M_0$ is a modification of its unique normal Stein filling. We call a partial filling of $M_0$ a complex manifold $X$ whose boundary components are $M_0$ as the pseudoconvex end and a SPCR-3-manifold $M_1$ as the pseudoconcave end. We show that the general problem of the embeddability of the pseudoconvex end of $X$ is equivalent 
to the rigidity of the possible partial fillings of $M_0$. At the same time, in Proposition 1 and also Theorem 3, we give methods to construct pseudoconcave surfaces that demonstrate that there is no rigidity in the partial fillings of an embeddable SPCR-3-manifold. \par
As a consequence, Theorem 1 of this paper answers negatively a question raised by Epstein and Henkin in [EpHe1-00] asking if  the embeddability of  SPCR-3-manifolds is preserved under complex cobordism. One reason for their  question arose from the study of the closedness for the $C^{\infty}$ topology of the set of small embeddable deformations a given CR-structure on a 3-manifold. In [EpHe2-00] one is led to analyse the embeddability of SPCR-3-manifolds $\partial X_-$ that are the boundaries  of  pseudoconcave surfaces $X_-$ which have inside another pseudoconcave surface $Y_-$ whose boundary $\partial Y_-$ is embeddable. It was thought that this might be enough to guarantee the embeddability of $\partial X_-$, one of the consequences of our results is that it is not. \par

As we shall see, one of the obstructions to get the embeddability of $\partial X_-$ is that the map giving the embedding $Y_-$ in general does not extend to a map of $X_-$. This has to do with the non-extendability of CR-functions on $\partial Y_-$ to holomorphic functions on $\mbox {Int}(X_- \setminus Y_-)$. We show this is possible using the method to construct pseudoconcave surfaces described by Proposition 1. This is a new phenomenon that is not possible in higher dimensions. On the other hand, Theorem 2 shows that the extendability of all CR-functions of $\partial Y_-$ and therefore the extendability of the map giving  the embedding of $Y_-$ to $X_-$ is not equivalent to the embeddability of $\partial X_-$ (though it is necessary by a result of Kohn and Rossi [KoRo65]).\par

Another interpretation for Theorem 2 is that complex cobordisms not preserving embeddability can be similar to complex cobordisms preserving embeddability. Our last result renforces that claim. Theorem 3 shows that in a family of complex-cobordism manifolds the property that the complex-cobordism preserves embeddability is not locally stable. This theorem also implies a result  on the stability of embeddability for CR-structures on a 3-dimensional manifold that are complex cobordant to a given embeddable SPCR-3-manifold.

\

{\bf Acknowledments} I would like to thank Professor Charles Epstein for generously sharing his knowledge with me and for his enthusiasm and Professor Siu for always inspiring and encouraging me.

\

\

\noindent {\bf SECTION 1}

\

\

Let $X$ be a compact complex (n+1)-manifold with boundary, such that the boundary components $M_i$ are all SPCR-(2n+1)-manifolds. Each of the boundary components of $X$, $M_i$, has an open neighborhood $D_i \subset X$ and a smooth function $\phi_i$ defined on $D_i$ such that $M_i=\{\phi_i=0\}$ and $\phi_i$ is strictly plurisubharmonic on $D_i \setminus M_i$. Moreover, $M_i$ will be one of the following: pseudoconvex end (resp. pseudoconcave end) if $\phi_i\le 0$ (resp.  $\phi_i\ge 0$) in a neighborhood of $M_i$.\par

\

{\bf Definition} A compact complex manifold $X$ as described above is called a (1,1)-convex-concave manifold. A (1,1)-convex-concave complex surface with no pseudoconvex (pseudoconcave) ends is called a pseudoconcave (resp. pseudoconvex) surface.

\

We analyze how the embeddability of one boundary component of a (1,1)-convex-concave manifold influences the embeddability of the other components. We introduce the notation of [EpHe1-00]: 

\

{\bf Definition} Two CR-manifolds $M_0$ and $M_1$ are called complex-cobordant if there is a (1,1)-convex-concave manifold $X$, the complex-cobordism manifold,  such that $\partial X=M_0 \cup M_1$ (as CR-manifolds), $M_0$ and $M_1$ are respectively the pseudoconvex and  pseudoconcave ends. If in addition there is a strictly plurisubharmonic function $\rho$ on $X$ with $M_0$ and $M_1$ as level sets then $M_0$ and $M_1$ are called strictly CR-cobordant.

\

Let $M_0$ and $M_1$ be two complex cobordant SPCR-3-manifolds with complex cobordism manifold $X$. What is the relation between the embeddability of $M_0$ and $M_1$? The first observation is that the embeddability of the pseudoconcave end $M_1$ implies the embeddability of the pseudoconvex end $M_0$ since in this case $M_0$ would be fillable. The relevant problem is the implication in the other direction, does the embeddability of the pseudoconvex end imply the embeddability of the pseudoconcave end? \par

A result of Epstein and Henkin [EpHe3-00] shows that if the SPCR-manifolds $M_0$ and $M_1$ are  strictly CR-cobordant then $M_0$ is embeddable iff $M_1$ is embeddable. This result also answers the local version of the question raised above, for in the (1,1)-convex-concave manifold $X$ there is a sufficiently small (1,1)-convex-concave collar of $M_0$ which is a strict CR-cobordism between its boundary components.\par

The embeddability of the the pseudoconvex end $M_0$ implying the embeddability of the pseudoconcave end $M_1$ has geometrical interpretations that motivated the question above as well as its answer. 

\

{\bf Definition} Let $M_0$ be an embeddable be SPCR-(2n+1)-manifold. We call a (1,1)-convex-concave manifold $X$ with $M_0$ as its only pseudoconvex end a partial filling of $M_0$. We say that we have rigidity on the partial fillings of $M_0$ if they are all modifications of subsets of the unique Stein normal space filling $M_0$. 

\

It follows from [Ro65] that for $n\ge 2$ all partial fillings of an SPCR-(2n+1)-manifold $M_0$ are rigid. The following describes the relationship between the rigidity of the possible partial fillings of $M_0$ and the embeddability of the pseudoconcave ends of those partial fillings. The embeddability of $M_0$ implies that $M_0$ bounds a unique normal Stein surface $S$. Assume that all pseudoconcave ends $M_i$ are embeddable and $S_i$ are their respective Stein fillings. Let $S'=X \amalg^i_{M_i} S_i$ be the complex space obtained by gluing $X$ to the $S_i$ along the $M_i$. Since $S'$ is a strictly pseudoconvex filling of $M_0$ it follows that $X$ is an open subset of $S$ or of one of its modifications. The result of Epstein and Henkin states, regardless of the dimension,  that in any partial filling of $M$ there is a sufficiently small collar of $M_0$ in $X$ that is rigid. We will demonstrate that partial fillings of an embeddable SPCR-3-manifold are far from rigid and in particular some will not extend to an actual filling. This is another difference between SPCR-3-manifolds and SPCR-(2n+1)-manifolds for $n\ge2$. \par

Lempert in [Le95] showed that every embeddable SPCR-3-manifold $M$ can be presented as the separating hypersurface in a projective surface $W$. $M$ separates $W$ in two components: the pseudoconcave $W_-$ with a positive curve $C \subset W_-$ and the pseudoconvex $W_+$. This implies that any embeddable SPCR-3-manifold $M_0$ is the boundary of a pseudoconcave surface $W_-$ moreover this pseudoconcave surface is embedded in some projective surface $W$ (note that there is no rigidity for the pseudoconcave fillings of $M_0$). Therefore the partial fillings of $M_0$ can be seen as extensions of the pseudoconcave surface $W_-$ that are also pseudoconcave surfaces. This will be the approach taken in this paper to produce examples of complex cobordisms.\par

The next lemma is one of the essential ingredients in  the construction of the examples on this paper.

\

{\bf Lemma 1} If $X$ is a nonsingular complex surface containing the normal crossing divisor $C=C_1 \cup ...\cup C_k$, with $C_i$ nonsingular complete curves with $C_i^2>0$. Then there exists a neighborhood $W$ of $C$ and a smooth strictly plurisubharmonic function $\varphi:W\setminus C \to \mathbb R$ such that:

\

\noindent

i) For every sequence $(x_n)_{n \in \mathbb N}$, $x_n \in W\setminus C$ converging to $x \in C$, $\lim_{x_n \to x}\varphi(x_n)=+ \infty$.

\

\noindent

ii) For $c \gg 0$, $X_c=\{x \in W: \varphi(x)\ge c \mbox { or } x\in C \}$ is a pseudoconcave neighborhood of $C$.

\

{\bf Proof:} We prove the statement for $C=C_1 \cup C_2$ with $C_1$ and $C_2$ intersecting transversely at the single point $C_1 \cap C_2={x_0}$, the general statement follows from the arguments for this case. Before starting the proof recall:

\

\noindent 1)  A function $f:U \to \mathbb R$, $U \subset \mathbb C ^n$, is strictly plurisubharmonic if the Levi form, the quadratic form associated with the complex Hessian $[\partial \bar \partial f]$
:

$$Lf \equiv \sum_{\alpha,\beta}\frac{\partial ^2f}{\partial \overline z_{\beta}\partial z_{\alpha}}dz_{\alpha}\otimes d\overline z_{\beta}$$
is positive definite. 

\

\noindent 2) If h is a hermitian metric on a line bundle $L$. Then:

$$-\partial \overline {\partial} \log h =\sum R_{\alpha,\overline \beta}dz_{\alpha}\wedge  d\overline z_{\beta}$$

\noindent where $\sum R_{\alpha,\overline \beta}dz_{\alpha}\wedge d\overline z_{\beta}$ is the Ricci form for the hermitian connection.

\

\noindent 3) A line bundle $L$ on a compact Kahler manifold $M$ is positive/negative if there is a hermitian metric having a positive/negative definite  matrix of coefficients of the Ricci form at every point of $M$ (such metric is called  positive/negative). For curves a line bundle is positive if its degree is positive.

\

We proceed by constructing the function $\phi:W\setminus C_1\cup C_2 \to \mathbb R$ with the desired properties.
Let $(U_i)_{i \in \mathbb N}$ be a covering of $X$ consisting of relatively compact  open subsets of $X$. Denote by $f_i$ the defining equations for $C_1 \cup C_2$ in $U_i$. Choose this covering such that $x_0 \in U_0$ and $f_0=z_1z_2$ and  for all other $U_i$  either i) $f_i=1$ or ii) $f_i=z_1$. \par

Since $C^2_1,C^2_2>0$ then $\mbox {deg}(O(C_1 +C_2)|_{C_1}), \mbox {deg}(O(C_1 +C_2)|_{C_2} )>0$ and by 3) we can find positive metrics $g^1$, $g^2$ on respectively $O(C_1+C_2)|_{C_1}$ and $O(C_1 + C_2)|_{C_2}$. To define a metric $g$ on $O(C_1 +C_2)|_{C_1 \cup C_2}$ we scale $g^1$ such that $g^1$ and $g^2$ agree at $x=x_0$ and make $g|_{C_1}=g^1$ and $g|_{C_2}=g^2$.  Let
 $h$ be a metric on $O(C_1+C_2)$ such that $h=g$ on $O(C_1+C_2)|_{C_1\cup C_2}$, denote $h_i=h|_{U_i}$.\par
The desired function $\phi$ is constructed by gluing local functions of the following form:

$$\phi_i=\phi|_{U_i\setminus \{f_i=0\}}=-\log(h_i |f_i|^2)$$

Since $\partial \bar \partial \phi_i=\partial \bar \partial (-\log(h_i))$ where $f_i\ne 0$, the properties of the Levi form of $\phi_i$, $L(\phi_i)$, follow from the properties of $L(-\log(h_i))$, i.e properties of the curvature of the metric. In particular, if we can find a neighborhood $W$ of $C_1 \cup C_2$ where the metric $h$ is positive then $\phi$ is strictly plurisubharmonic on $W\setminus C_1 \cup C_2$. The next step then is to make the metric $h$ positive in a neighborhood of $C_1\cup C_2$. \par
First, we modify $h$ to make the metric positive near $x_0$. We achieve this by making the matrix  $[\partial \bar \partial \log h_0]$ negative definite at $x=x_0$, but preserving $h|_{C_1 \cup C_2}=g$. Since $h_0|_{\{z_1=0\}}=g^1|_{C_1\cap U_0}$ and $h_0|_{\{z_2=0\}}=g^1|_{C_2\cap U_0}$ and both metrics $g^1$, $g^2$ are positive, it follows that the diagonal elements of  $[\partial \bar \partial \log h_0](x_0)$ are negative.
We can assume that the covering is such that there is a ball $B_{x_0}(r)$ satisfying $B_{x_0}(r)\cap U_i=\emptyset$, $i\ne 0$. Change $h_0$ to $h_0e^{2Re(az_1\bar z_2).\rho}$ where 
$a=-\frac{\partial ^2 \log h_0 (x_0)}{\partial \bar z_2 \partial z_1}$ and $\rho$ a smooth function on $U_0$ with value $1$ on $B_{x_0}(\frac{r}{2})$ and value $0$ outside $B_{x_0}(r)$. The new $h_0$ has a negative definite diagonal complex Hessian since its entries are equal to the diagonal entries of the complex Hessian for the previous $h_0$.

Second, we modify $h$ to make the metric also positive away from $x_0$ but near the curve $C_1 \cup C_2$. This modification is done by changing the $h_i$ to:

$$h'_i=\frac {h_i}{1+ch|f_i|^2}$$

This type of modification of $h_i$ was done for the smooth case by [Sc73]. The complex Hessian of $\log h'_i$ is given by $[\partial \bar \partial \log h'_i]=[\partial \bar \partial \log h_i] -[\partial \bar \partial \log (1+ch|f_i|^2)]$. We show that by increasing $c$ we can make $[\partial \bar \partial \log h'_i]$ negative definite for all points $x \in C_1 \cup C_2$. The contribution of $[\partial \bar \partial \log (1+ch|f_i|^2)]$ is:

\

\noindent 1) $\partial \bar \partial \log (1+ch_i|z_1|^2)=ch_idz_1 \wedge d\bar z_1$ on $U_i\cap \{z_1=0\}$, $i \ne 0$.

\

\noindent 2) $\partial \bar \partial \log (1+ch_0|z_1z_2|^2)=ch_0|z_2|^2 dz_1 \wedge d\bar z_1 + ch_0|z_1|^2 dz_2 \wedge d\bar z_2$ on $U_0\cap \{z_1z_2=0\}$.

\

There is a ball $B_{x_0}(r)\subset U_0$ such that $[\partial \bar \partial \log h_0]$ is negative definite on $\{z_1z_2=0\} \cap B_{x_0}(r)$, since  $[\partial \bar \partial \log h_0](x_0)$ is negative definite. Note that it follows from 2) that for any $c>0$, the matrix  $[\partial \bar \partial \log h'_0]$ differs from $[\partial \bar \partial \log h'_0]$ only on  the diagonals, which become more negative (the negativity of the diagonal entries follow from the positivity of the metrics $g^1$ and $g^2$). So $[\partial \bar \partial \log h_0]$ is negative definite on $\{z_1z_2=0\} \cap B_{x_0}(r)$ for all $c>0$. \par

Without loss of generality we can assume that  any $x\in C_1 \cup C_2$ not in $B_{x_0}(r)$ is  in $U_i$, $i\ne 0$. We show that when $c$ sufficiently big there is a neighborhood $V\subset X$ of $x$ where $[\partial \bar \partial \log h'_i]$ is negative definite. Since the metrics $g^1$ and $g^2$ are positive, it follows that the diagonal entry  $\frac{\partial ^2 \log h_i}{\partial \bar z_2 \partial z_2}(x)$ of  $[\partial \bar \partial \log h_i](x)$ is negative. By 1) we see that $[\partial \bar \partial \log h'_i]$ differs from $[\partial \bar \partial \log h_i](x)$ only on the other diagonal which  can be made as negative as desired by increasing $c$. Hence there is a $c$ for which $[\partial \bar \partial \log h'_i](x)$ is negative definite. Then the  existence of a metric $h'$ which is positive in a neighborhood of $C_1 \cup C_2$ follows from a compactness argument. This establishes i) and ii) follows immediately.

\

The next proposition  shows that pseudoconcave surfaces can have two nonintersecting positive curves. An application of this result is done in the proof of Theorem 1. We note that the analogous result does not hold in higher dimensions (see the remark after Corollary 1)

\

{\bf Proposition 1} There are pseudoconcave surfaces with positive curves $C_i$, i.e $C^2_i>0$, that do not intersect.
Moreover one can assume that each of these curves has an embeddable pseudoconcave neighborhood.

\
 
{\bf Proof} We will construct a pseudoconcave surface $S$
containing two positive curves $C_1$ and $C_3$ that do not intersect, the proof of the general statement follows directly from the argument described below. \par

Let $C_1$, $C_2$ and $C_3$ be 3 positive curves
in the compact complex surfaces $X_1$, $X_2$ and $X_3$ respectively. First, we glue a neighborhood of $C_1$ in $X_1$  with a neighborhood of $C_2$ in $X_2$ obtaining a surface $X_4$ containing $C_1$ and $C_2$ intersecting transversely at one point. Pick a  point $p_1 \in C_1$. Let $U_1$ be a subset of
$X_1$ containing $p_1$ and $\phi$ a biholomorphism $\phi: U_1 \to \Delta \times
\Delta$, $\Delta$ the disc of radius 1, such that $\phi(C_1 \cap U_1) =0 \times \Delta$. Let $U_2$ be an open subset
of $X_1$ such that $C_1 \subset U_1 \cup U_2$ and $\phi(U_1\cap U_2) \subset
\Delta_\frac {1}{2}\times \Delta$, $\Delta_\frac {1}{2}$ the disc of radius 1/2.  \par

Do the same for $C_2$ in $X_2$. Pick a point $p_2 \in C_2$, an open subset $V_1 \subset X_2$ and $\varphi$ a biholomorphism $\varphi: V_1 \to \Delta \times \Delta$ such that $\varphi(C_2 \cap V_1)=\Delta \times 0$. Pick $V_2$ such that $C_2 \subset V_1 \cup V_2$ and $\varphi(V_1 \cap V_2) \subset \Delta \times \Delta_\frac{1}{2}$.

The surface $X_4$ is the surface given by the open sets $\Delta \times \Delta$, $U_2$ and
$V_2$ with the gluings of $U_2$ and $V_2$ with $\Delta \times \Delta$ given by respectively $\phi|_{U_1 \cap U_2}$ and $\varphi|_{V_1 \cap V_2}$ .\par

Now repeat the same argument for the curve $C_1 \cup C_2$ in $X_4$ and the curve $C_3$ in $X_3$, picking a point  in $C_2$ that is not the crossing point. The resulting surface $X$ has two positive curves, $C_1$ and $C_3$ not intersecting and the neighborhood germs of $C_1$, $C_2$ and $C_3$ in $X$ will be the initial neighborhood germs in $X_1$, $X_2$ an $X_3$ respectively. The surface $X$ and the curves $C_1$, $C_2$ and $C_3$ satisfy the conditions of Lemma 1 hence there is a pseudoconcave neighborhood of $C_1 \cup C_2 \cup C_3$ on $X$.

\

{\bf Theorem 1} Embeddability of CR-3-manifolds is not a complex-cobordism invariant.

\

{\bf Proof:} Let $M_0$ and $M_1$ be two complex-cobordant SPCR-3-manifolds and $X$ the complex-cobordism manifold with $M_0$, $M_1$ as the pseudoconvex and pseudoconcave ends respectively. The pseudoconcave end  $M_1$ being embeddable is the same as $M_1$ being fillable. This in turn, implies that $X$ is a subset of a modification of the unique Stein normal filling of $M_0$. Kohn and Rossi [KoRo65] showed that if a SPCR-manifold $M$ bounds a Stein normal space $S$, then any CR-function on $M$ extends to an holomorphic function of $S$. Therefore the embeddability of $M_1$ implies that any CR-function on $M_0$ must extend to an holomorphic function on the interior of $X$.\par
The surface constructed in the proof of proposition 1 contains  a pseudoconcave neighborhood $W$ of the connected chain of positive curves $C_1\cup C_2 \cup C_3$. Recall that by the construction of the surface, $C_1$ has an embeddable pseudoconcave neighborhood $W_1 \subset W$. The complex manifold $X=W \setminus W_1$ has the embeddable pseudoconvex boundary component $M_0=\partial W_1$, and the pseudoconcave boundary component $M_1=\partial W$. Since $M_0$ is embeddable it must have non constant CR-functions. On the other hand, $X$ contains a pseudoconcave neighborhood of the positive curve $C_3$ in its interior, therefore all holomorphic functions on the $\mbox {Int}(X)$ must be constant . This follows from the maximum principle and the fact that the pseudoconcavity of $W_3$ implies that the set $\hat W_3=\{x\in X||f(x)|\le \sup_{y \in W_3}|f(y)| \mbox { for all holomorphic functions on Int(X)}\}$ satisfies $\overline W_3 \subset \mbox {Int}(\hat W_3)$ [AnGr62].  The constancy of all holomorphic functions on $\mbox {Int}(X)$ implies that not all CR-functions on $M_0$ extend to holomorphic functions on the interior of $X$, hence by the previous paragraph $M_1$ is not embeddable.

\

{\bf Corollary 1} (of Proposition 1) There are (1,1)-convex-concave manifolds of complex dimension 2 with more than one pseudoconvex end.

\
\newpage

{\bf Remark} The phenomenon described in Corollary 1 is only possible in dimension 2 since Rossi in [Ro65] proved that a (1,1)-convex-concave manifold of complex dimension $n$, $n>2$, has at most one pseudoconvex end. The result analogous to Proposition 1 would propose, for all dimensions, the existence of a pseudoconcave complex manifold $Y$ containg two complete hypersurfaces having disjoint pseudoconcave neighborhoods. For the sake of completness, we give a quick argument that shows that the generalized proposition does not hold if $\mbox {dim }Y>2$, a fact that also follows imediately from Rossi's result. Let $X$ be the compact complex space obtained by gluing to $Y$ the filling of $\partial Y$, which exists since $\mbox {dim }Y>2$. $X$ has two pseudoconcave neighborhoods $W_1$ and $W_2$, of respectively the complete hypersurfaces $H_1$ and $H_2$, satisfying $W_1 \cap W_2=\emptyset$. Then $X\setminus \overline W_1$ is a strictly pseudoconvex complex space containing the pseudoconcave neighborhood $W_2$ of $H_2$, which is an impossibility. On one hand $X\setminus \overline W_1$ has non-constant holomorphic functions since it is an open holomorphic convex space, and on the other hand it has no non-constant holomorphic functions since it contains the pseudoconcave neighborhood $W_2$.

\

\

\noindent {\bf SECTION 2}

\

\

The examples of complex cobordisms described in the proof of Theorem 1 are very degenerate with respect to preserving embeddability. In particular, the complex-cobordism manifold $X$ has positive curves inside and therefore no non-constant CR-function on the pseudoconvex end $M_0$  extends to a holomorphic function on $\mbox {Int}(X)$. We proceed to show that complex cobordisms not preserving embeddability can be quite similar to the ones that do. \par

The most natural candidate for a great divide between complex-cobordisms preserving and not preserving embeddability would be on the extendability of the CR-functions of the pseudoconvex end. The next result shows this divide does not occur.

\

{\bf Theorem 2} Let $M_0$ be an embeddable SPCR-3-manifold and $M_1$ be complex-cobordant to $M_0$ with complex-cobordism manifold $X$ having $M_0$ as the pseudoconvex end. The extendability of the CR-functions on $M_0$ to holomorphic functions on $\mbox {Int}(X)$ does not imply that $M_1$ is embeddable.

\

{\bf Proof} As was explained in the proof of theorem 1, the extendability of the CR-functions is a necessary condition for the embeddability of $M_1$, we show that it is not sufficient.\par
To prove the theorem we proceed to construct a pseudoconcave surface containing three $\mathbb P^1$ where two of the $\mathbb P^1$ do not intersect but each of these two intersects the third $\mathbb P^1$ at one point (which is not the same). We also ask for the germs of the neighborhoods to be equivalent to the neighborhood of the linear $\mathbb P^1 \subset \mathbb P^2$ (we denote such curves as the linear $\mathbb P^1$'s). One way to construct this surface is to use the method described in proposition 1. We now give another way to do this construction. \par

Consider three linear $\mathbb P^1 \subset \mathbb P^2$ $C_1$, $C_2$ and $C_3$ with $U_1$, $U_2$ and $U_3$, tubular neighborhoods of each of these curves such that $U_1\cap U_2 \cap U_3=\emptyset$.  Let $W'$ be the complex surface that results from $U_1 \cup U_2 \cup U_3 \subset \mathbb P^2$ by separating $U_2$ from $U_3$. The curves in $W'$ coming from the three linear $\mathbb P^1 \subset \mathbb P^2$ and the surface $W'$ satisfy the conditions stated in the above paragraph. Denote the three linear $\mathbb P^1$'s in $W'$ also  by $C_1$, $C_2$ and $C_3$ ($C_2 \cap C_3=\emptyset$). Lemma 1 implies that there is a pseudoconcave surface $W\subset W'$ containing the curve $C_1 \cup C_2 \cup C_3$. Since $W$ contains 2 disjoint positive curves $C_2$ and $C_3$, $W$ is not an embeddable pseudoconcave surface , which in turn implies that $\partial W$ is also not embeddable.\par

After making the identification of $\mathbb P^2 \setminus C_1$ with $\mathbb C^2$, consider a sufficiently large ball such that its complement $W_1$ in $\mathbb P^2$ is contained in $U_1 \cap W$ ($W$ lives in $W'$ but $U_1 \cap W$ is a subset of both $\mathbb P^2$ and $W'$). Clearly all CR-functions of the boundary of the ball extend to the interior of the ball and therefore to $U_1 \cup U_2 \cup U_3 \setminus W_1$. If one considers $W_1$ as a pseudoconcave surface in $W$,  this in turn implies that all CR-functions on $\partial W_1$ extend to holomorphic functions on $\mbox{Int}(W\setminus W_1)$ and the result follows.

\

In order to express a consequence of Theorem 2, we introduced the two weakened versions of embeddability of a pseudoconcave surface $X_-$ described in [EpHe2-00].

\

{\bf Definition:} Let $X_-$ be a pseudoconcave surface with a positive curve $Z$, $(X_-, Z)$ is called weakly embeddable if the holomorphic map defined by the sections of $H^0(X_-,O(dZ))$, $\varphi_d:X_- \to \mathbb P^n$, is injective in some neighborhood of $Z$ for sufficiently large $d$. The pair $(X_-, Z)$ is called almost embeddable if in addition the map $\varphi_d:X_- \to \mathbb P^n$ is an embedding outside a proper analytic subset of $X_-\setminus Z$..

\

In [EpHe2-00] the authors ask about the relationship between these two notions of weakened embeddability. The interest sprouted from the fact that one has an operative numerical criterion involving the dimensions $h^0(X_-,O(dZ))$ characterizing weak embeddability and that almost embeddability of $X_-$ implies embeddability of $\partial X_-$. In the proof of Theorem 2 we constructed the weakly pseudoconcave surface $(W,C_1)$ such that for any $d$ the maps $\varphi_d:W \to \mathbb P^n$ can not distinguish any pair points in $W$ that came from the same point in $U_2 \cap U_3$ and we have: 

\

{\bf Corollary 2}  The weak embeddability of a pseudoconcave surface $(X_-,Z)$ does not imply that $(X_-,Z)$  is almost embeddable. 

\

\

\noindent {\bf SECTION 3}

\

\

As discussed above, there are geometrical consequences of the extendability of CR-functions from the pseudoconcave end. As a particular example of this, we have the next proposition that is used the proof of Theorem 3.

\

{\bf Proposition 2} Let $X$ be a pseudoconcave surface having two positive nonsingular curves $C_1$ and $C_2$ with one normal crossing. Assume that $C_1=\mathbb P^1$, $C_1^2=1$ and $C_1$ has an embeddable pseudoconcave neighborhood $X_1\subset X$. If all CR-functions on $\partial X_1$ extend to holomorphic functions on $\mbox{Int}(X\setminus X_1)$ then $C_2=\mathbb P^1$ and the germs of the neighborhoods of $C_1$ and $C_2$ in $X$ are equivalent to the germ of the neighborhood of the linear $\mathbb P^1\subset \mathbb P^2$.

\

{\bf Proof} By  [MoRo81] we know that the neighborhood germ of $C_1$ in $X$ is equivalent to neighborhood germ of the linear $\mathbb P^1 \subset \mathbb P^2$. For convenience of the rest of our argument, we give a proof of this fact. By hypothesis $X_1$ is embeddable, hence $X_1$ is contained in a smooth projective surface $S$. $S$ must be a rational surface since it contains the rational curve $C_1$ with positive self-intersection. 
From the exact sequence $0 \to O_S \to O_S(C_1) \to O_{C_1}(C_1)\to 0$ and the vanishing $H^1(S,O_S)=0$, it follows that $h^0(S,O_S(C_1))=3$. Any basis of   $H^0(S,O_S(C_1))$ defines a map $f:S \to \mathbb P^2$ which by the Stein factorization is a composition of a  map with connected fibers and a finite map. The connected map contracts exactly the curves not intersecting $C_1$. The finite map is the identity map since the preimage a general point in $\mathbb P^2$ is one point. The preimage of a general point in $\mathbb P^2$ is given by the intersection of two smooth sections  in  $H^0(S,O_S(C_1))$, which is one point since $C^2_1=1$. From the previous, it follows that the map $f:S \to \mathbb P^2$ is an embedding of a neighborhood of $C_1$ and the image of $C_1$ by $f$ is a linear $\mathbb P^1 \subset \mathbb P^2$.\par

The map $f: S \to \mathbb P^2$ is given by $g(x)=[s_0(x):s_1(x):s_2(x)]$, where $s_0$, $s_1$ and $s_2$ are a basis of $H^0(S,O(C_1))\subset H^0(X_1,O(C_1))$. We are interested in extending $f|_{X_1}$ to a map $f':\mbox {Int} (X) \to \mathbb P^2$. To achieve this goal, we show that the sections $s_0,s_1,s_2\in H^0(X_1,O(C_1))$ can be extended to sections $s'_0,s'_1,s'_2\in H^0(\mbox {Int}(X),O(C_1))$.  \par

A set of defining equations for $C_1$ in $\mbox {Int}(X)$ give a section $s'_0\in  H^0(\mbox {Int}(X),O(C_1))$. Without loss of generality we can assume that $s_0=s'_0|_{X_1}$. To extend the other sections $s_1$ and $s_2$, we use the equalities $s_i=\frac {s_i}{s_0}s_0$, $i=1,2$. Note that $\frac {s_1}{s_0}$ and $\frac{s_2}{s_0}$ are holomorphic functions on $S \setminus C_1$. The restrictions of these functions, $\frac {s_1}{s_0}|_{_{\partial X_1}}$ and $\frac {s_2}{s_0}|_{_{\partial X_1}}$, are CR-functions on $\partial X_1$. The assumptions on $X_1$ imply that the CR-functions $\frac {s_1}{s_0}|_{_{\partial X_1}}$ and $\frac {s_2}{s_0}|_{_{\partial X_1}}$ extend. Therefore we obtain two holomorphic functions on  $\mbox {Int} (X \setminus C_1)$ denoted by $\frac {s'_1}{s'_0}$, $\frac{s'_2}{s'_0}$ that extend $\frac {s_1}{s_0}|_{X_1\setminus C_1}$ and $\frac{s_2}{s_0}|_{X_1\setminus C_1}$. The extensions $s'_1$ and $s'_2$ follow from the extended identities $s'_i=\frac {s'_i}{s'_0}s'_0$, $i=1,2$.\par

In conclusion, $h^0(\mbox {Int}(X), O(C_1))\ge 3$ and its sections $s'_0$, $s'_1$, $s'_2$ give an extension of the map $f|_{X_1}:X_1 \to \mathbb P^2$ to a map $f':\mbox {Int}(X) \to \mathbb P^2$. As an immediate consequence we have that $C_2$ must be a rational curve. The map $f'|_{C_2}:C_2 \to \mathbb P^2$ is given $s'_0|_{C_2}, s'_1|_{C_2},s'_2|_{C_2} \in H^0(C_2,O(C_1)|_{C_2})$ and since the image $f'(C_2)$ is not a point, it follows that $h^0(C_2,O(C_1)|_{C_2})\ge 2$. $C_2$ must be a rational curve, since only rational curves have line bundles of degree 1 with two linearly independent sections.\par

We want to show that $f'$ embeds a neighborhood $X_2 \subset X$ of $C_2$. As a first step, we get that $f'|_{C_2}: C_2 \to \mathbb P^2$ is an embedding and its image is a linear $\mathbb P^1$, since $\deg f'(C_2).\deg f'=\deg_{C_2} {f'}^*O(1)=C_1.C_2=1$. As a second step, we obtain that the differential $df'_*$ has full rank at all points of $C_2$. To see this, we use the commutative diagram:

\begin{equation}\label{a}
\xymatrix{0 \ar [r]  &T_{C_2} \ar [r]\ar [d]^{id} & T_{X|_{C_2}} \ar [r]\ar [d]^{df'_*}&N_{C_2}\ar [r]\ar [d]^{\overline {df'_*}}&0\\
         0\ar [r]&{f'}^*T_{f'(C_2)}\ar [r]&{f'}^*T_{\mathbb P^2|_{f'(C_2)}}\ar [r]&{f'}^*N_{f'(C_2)}\ar [r]&0}
\end{equation}

After the isomorphisms $N_{C_2}\simeq O(n)$ and ${f'}^*N_{f'(C_2)} \simeq O(1)$, morphism $\overline {df'_*}$ of (1) becomes $\overline {df'_*}: O(n) \to O(1)$ and hence  $\overline {df_*}$ is a section of $H^0(\mathbb P^1,O(1-n))$. Since $H^0(\mathbb P^1,O(1-n))$ is trivial for $n>1$ and $\mathbb C$ for $n=1$, it follows that for $n>1$ the morphism $\overline df'_*$ must be the trivial morphism and  for $n=1$ the morphism $\overline df'_*$ can be trivial or an isomorphism. The morphism $\overline {df'_*}$ can not be trivial because in a sufficiently neighborhood $U \subset X_1$ of $p_0=C_1 \cap C_2$,  $f'|_{U}=f|_U$ and $f|_U$ is an embedding. Hence $n=1$, i.e $C^2_2=1$ and $df'_*$ has full rank at all points of $C_2$. The implicit function theorem implies the desired embedding of a neighborhood of $C_2$.

\

The last result has two interpretations. One is that the property that a complex manifold gives a complex-cobordism preserving embeddability is not stable under small deformations of the complex-cobordism manifold. The other is related to the understanding of  the structure of the set of embeddable SPCR-3-manifolds which is one of the main goals of the subject. The result implies that the behavior of embeddability for  SPCR-3-manifolds complex cobordant to a given embeddable SPCR-3- manifold $M_0$ share some of the same pathologies as in the general case, where there are no restrictions.

\

{\bf Theorem 3} In a family of complex-cobordism manifolds whose pseudoconvex ends are embeddable the property that the respective complex-cobordisms preserve embeddability is not locally stable.

\

{\bf Proof} One proceeds by constructing a counterexample to the stability of the preservation of embeddability for complex-cobordisms. \par

Let $C_1$ and $C_2$ be two distinct linear $\mathbb {P}^1$'s in $\mathbb P^2$ and $x_0 =C_1 \cap C_2$. Let $U_0$, $U_1$ and $U_2$ be a covering of a neighborhood $U$ of $C=C_1\cup C_2$ such that  $x_0 \in U_0$, $C_1 \subset U_0 \cup U_1$, $C_2 \subset U_0 \cup U_2$, and $U_1 \cap U_2=\emptyset$. We can choose  $U_0\simeq \Delta_1 \times \Delta_1$, $U_2\simeq\Delta_1 \times \Delta_1$ satisfying $U_0 \cap C_2=\Delta_1 \times 0$, $U_2 \cap C_2=\Delta_1 \times 0$,  $U_0 \cap U_1\subset \Delta_{\frac{1}{2}}\times \Delta_1 \subset U_0$ and $U_0 \cap U_2\subset \{z\in \Delta_1||z|>\frac{1}{2}\} \times \Delta_1 \subset U_0$ where $\Delta_r$ is a disc with radius r. Moreover, under the identifications of $U_0$ and $U_2$ with $\Delta_1 \times \Delta_1$ done above, the gluing of $U_0$ with $U_2$ is given by a fiberwise map f, with respect to the fibers of the projection to the first factor. [MoRo81] shows that there exist holomorphic families of gluings of $U_0$ and $U_2$ described by fiberwise maps parametrized by $t$ such that: for $t=0$ one has the gluing induced by $f$ ; for $t \ne 0$ the gluings give rise to open surfaces containing $C_2$ with the same normal bundle but the embedding of $C_2$ is inequivalent to the the embedding of $C_2$ for $t=0$ (see the remark below). As noted above, these non-standard germs are not filable.\par

Let $\omega:\cal V \to \Delta$ be the family of surfaces , with ${\cal V}_t =V_t$ and $V_0=U$, obtained by keeping the gluing of $U_0$ with $U_1$  but changing the gluing of $U_0$ with $U_2$  as  in the previous paragraph. Each member $V_t$ of the family has the curve $C$ embedded with the same normal bundle. Hence the tubular neighborhood of $C$ in all the $V_t$ is diffeomorphic to the tubular neighborhood $T$ of $C$ in $V_0$. There is a smooth map $\phi:T \times \Delta \to \cal V$, with each $\phi_t:T \to T_t=\phi(T\times t)\subset V_t$ a diffeomorphism from $T$ to a tubular neighborhood of $C$ in $V_t$. The family of surfaces $T_t$ can be described by the variation of the complex structure on $T$, induced by the diffeomorphisms $\phi_t$.\par

By Lemma 1 one can construct a strictly plurisubharmonic function $g:V_0 \setminus C \to \mathbb R$ satisfying $\{ x\in V_0 : g(x)\ge c\} \subset T$ for $c \gg 0$. Denote $S_c=\{ x\in T : g(x)= c\}$. Fix some $c \gg 0$, after possibly shrinking $\Delta$, one can assume that for all $t\in \Delta$, $g_t=g\circ \phi_t^{-1}:T_t\setminus C \to \mathbb R$ is strictly plurisubharmonic on a neighborhood of $\phi_t(S_c)$. Let $W_t\subset V_t$ be the set coming from the union of $U_0$ and $U_1$, by construction $W_t=W_0$ for all $t$. Repeating the arguments just described, one can pick a SPCR-3-manifold $M_0 \subset W_0 \cap \{ x\in T : g(x)> c\}$ bounding a neighborhood of $C_1$ such that $\phi_t(M_0)=M_{0t}\subset W_t$ is also an SPCR-3-manifold.\par

In conclusion, one has a family of complex-cobordism manifolds $\omega':\cal X \to \Delta$ such that each member $X_t$ is the (1,1)-convex-concave submanifold of $V_t$ with pseudoconvex end $M_{0t}$ and pseudoconcave end $M_{1t}$. By construction $M_{0t}$ is embeddable but  $M_{1t}$, for $t\ne 0$  is not embeddable because if it were the previous proposition would imply that the neighborhood of $C_2$ would be equivalent to the neighborhood germ of the linear $\mathbb P^1 \subset \mathbb P^2$. \par

Note that this proof shows that the property of the extendability of the CR-functions of the pseudoconvex end to holomorphic functions on the complex-cobordism manifold is also not locally stable in families of complex-cobordisms.

\

{\bf Remark} The paper [MoRo81] studies the space of nonequivalent holomorphic neighborhood retracts (HNR) of a curve $C$ with normal bundle $L$. An embedding of the  curve $C$ is a HNR if the neighborhood germ of $C$ has a holomorphic fibration over $C$. If $U$ and $V$ are a covering of $C$ such that $U$ is a disc and $U\cap V$ is an annulus, then we can describe all HNR of $C$ with normal bundle $L$, using the normal coordinates $u$, $v$ in respectively $U$, $V$, by:

$$ v=f(z,u)=g(z)(u+\sum_{i\ge 1}f_i(z)u^{i+1})$$

Where $g,f_i\in O(U\cap V)$ and $g$ is the transition function for the normal bundle $L$. To the HNR described above we can associate a point in $E(\bigoplus_{i>0}H^1(C,L^{-i}))$, the weighted projective space with the $\mathbb C^*$ action $t.(x_1,x_2,,...,x_n,...)=(tx_1,t^2x_2,...,t^nx_n,...)$. This is done by considering the $f_i$ as cocycles representing classes of $H^1(C,L^{-i})$. [MoRo81] shows that the neighborhood germs of two HNRs of $C$ with normal bundle $L$ are equivalent only if the corresponding points in  $E(\bigoplus_{i>0}H^1(C,L^{-i}))$ coincide. Since $H^1(\mathbb P^1,O(-i))\ne0$ for $i\ge 2$ the moduli space for HNRs is infinite dimensional.

\

{\bf Corollary 3} Let $M_0$ be an embeddable SPCR-3-manifold. The embeddability of SPCR-3-manifolds complex-cobordant to $M_0$ is not stable, for small deformations of the CR-structure preserving the property of being complex-cobordant to $M_0$.

\

{\bf Proof} In the proof of theorem 3 all the SPCR-3-manifolds $M_{1t}$ are complex cobordant to $S^3$ with the canonical embeddable CR-structure.

\

{\bf Remark} It would be very interesting to construct an example of a complex cobordidm manifold $X$ with $\partial X= M_0 \amalg M_1$, such that $M_0\simeq M_1\simeq S^3$, $M_0$ is fillable but $M_1$ is not.

\bigskip

\noindent B. de Oliveira, Harvard University and University of Pennsylvania. 

\

\noindent Harvard University,
 
\noindent Department of Mathematics, 

\noindent One Oxford Street, 

\noindent Cambridge, MA 02138

\

\noindent bdeolive@math.harvard.edu, bdeolive@math.upenn.edu

\bigskip


\begin{thebibliography}{60}


\bibitem [AnGr62]{AnGr62}{\bf Andreotti, A., Grauert, H.} {\em Theoremes de finitude pour la cohomologie des espaces complexes.\/ } Bull. Soc. Math. France, 90,(1962),193--259.
 
\bibitem [Bl94]{Bl94}{\bf Bland, J.} {\em Contact geometry and CR structures on $S^3$.\/ } Acta Math., 172,(1994),1--49.


\bibitem[Bou74]{Bou74}{\bf Boutet de Monvel } {\em Integration des equations de Cauchy-Riemann induites formelles.\/ }  Seminaire Goulaoic-Lions-Schwartz, Expose IX, (1974)

\bibitem[BuEp90]{BuEp90}{\bf Burns, D., Epstein, E. } 
{\em Embeddability for three-dimensional CR-manifolds.\/ } Jour. Am. Math. Soc, 3,(1990),809--841.

\bibitem[EpHe1-00]{EpHe1} {\bf Epstein, C., Henkin, G. } 
{\em Embeddings for 3-dimensional CR-manifolds.\/ }  Progress in Math., 188, Birkhauser (2000).

\bibitem[EpHe2-00]{EpHe2}{\bf Epstein, C., Henkin, G. } {\em Stability of embeddings for pseudoconcave surfaces and their boundaries.\/ }  Acta Math., 185:2, (2000), 161-237.

\bibitem[EpHe3-00]{EpHe3}{\bf Epstein, C., Henkin, G. } {\em Can a good manifold come to a bad end?.\/ }  Preprint 2000.

\bibitem[HaLa75]{HaLa75}{\bf Harvey, R., Lawson, H. } {\em On the boundaries of complex analytic varieties I.\/ } Ann. Math., 102,(1975),223--295.

\bibitem[Ha69]{Ha69}{\bf Hartshorne , R.} {\em Curves with high self-intersection on algebraic surfaces\/ } I.H.E.S, 36,(1969),111--125.

\bibitem [Ko63] {Ko63} {\bf Kohn, J.} {\em The range of the Cauchy-Riemann operators. \/ } Duke Math. J. , 53,(1986),525-545.


\bibitem[KoRo65]{KoRo65}{\bf Kohn, J., Rossi, H.}{\em On the extension of holomorphic functions from the boundary of a complex manifold. \/ } Ann. Math., 81,(1965),451-472.


\bibitem[Le94]{Le94}{\bf Lempert, L. } {\em Embeddings of three dimensional Cauchy-Riemann manifolds.\/ } Math. Ann., 300,(1994),1--15.


\bibitem[Le95]{Le95}{\bf Lempert, L. } {\em Algebraic approximations in analytic geometry.\/ } Inv. Math, 121,(1995),335--354.

\bibitem[Ro65]{Ro65}{\bf Rossi, H. } {\em Attaching analytic spaces to an analytic space along a pseudoconcave boundary.\/ } "Proc. Conf. on Complex Manifolds" (A.Aeppli, E. Calabi, Rohrl eds.)", Springer-Verlag,New York, Berlin, Heidelberg, 1965.
 


\bibitem[MoRo81]{M-R81}{\bf Morrow, J., Rossi, H. } {\em Some general results on equivalence of embeddings.\/ }  Princeton Annals, Proceedings of 1979 Princeton Conference on Complex Analysis, 1981, 299--325.

\bibitem[Sc73]{Sc73}{\bf Schneider, M. } {\em Uber eine Vermutung von Hartshorne.\/ } Math. Ann, 201,(1973),221--229.

\end{thebibliography}
\end{document}